\documentclass[12pt]{amsart}
\usepackage{amssymb,latexsym}
\newdimen\AAdi%
\newbox\AAbo%
\def\AAk#1#2{\s_etbox\AAbo=\hbox{#2}\AAdi=\wd\AAbo\kern#1\AAdi{}}%
\def\AAr#1#2#3{\s_etbox\AAbo=\hbox{#2}\AAdi=\ht\AAbo\raise#1\AAdi\hbox{#3}}%
\font\tenmsb=msbm10 at 12pt \font\sevenmsb=msbm7 at 8pt
\font\fivemsb=msbm5 at 6pt
\newfam\msbfam
\textfont\msbfam=\tenmsb \scriptfont\msbfam=\sevenmsb
\scriptscriptfont\msbfam=\fivemsb
\def\Bbb#1{{\tenmsb\fam\msbfam#1}}

\begin{document}

\newcommand{\D}{D \hskip -2.8mm \slash}
\newcommand{\pp}{\overline{\psi}}
\newcommand{\e}{\overline{\varepsilon}}
\newcommand{\hh}{\sqrt{h}d^2x}
\newcommand{\ii}{\frac{1}{2}\int_M}
\renewcommand{\theequation}{\thesection.\arabic{equation}}
\newcommand{\wb}{\widetilde{\nabla}_{e_\beta}}
\newcommand{\wa}{\widetilde{\nabla}_{e_\alpha}}

\newtheorem{thm}{Theorem}
\newtheorem{lem}{Lemma}
\newtheorem{cor}{Corollary}
\newtheorem{rem}{Remark}
\newtheorem{pro}{Proposition}
\newtheorem{defi}{Definition}
\newcommand{\noi}{\noindent}
\newcommand{\dis}{\displaystyle}
\newcommand{\mint}{-\!\!\!\!\!\!\int}
\newcommand{\ba}{\begin{array}}
\newcommand{\ea}{\end{array}}
\def \tf{\tilde{f}}
\def\cqfd{%
\mbox{ }%
\nolinebreak%
\hfill%
\rule{2mm} {2mm}%
\medbreak%
\par%
}
\def \pr {\noindent {\it Proof.} }
\def \rmk {\noindent {\it Remark} }
\def \esp {\hspace{4mm}}
\def \dsp {\hspace{2mm}}
\def \ssp {\hspace{1mm}}
\def\n{\nabla}
\def\RR{\Bbb R}\def\R{\Bbb R}
\def\C{\Bbb C}
\def\B{\Bbb B}
\def\N{\Bbb N}
\def\Q{\Bbb Q}
\def\Z{\Bbb Z}
\def\EE{\Bbb E}
\def\H{\Bbb H}
\def\SS{\Bbb S}\def\S{\Bbb S}
\def \c {{\bf C}}
\def \Z {{\bf Z}}
\def \Q {{\bf Q}}
\def \a {\alpha}
\def \b {\beta}
\def \d {\delta}
\def \e {\epsilon}
\def \G {\Gamma}
\def \g {\gamma}
\def \l {\lambda}
\def \L {\Lambda}
\def \O {\Omega}
\def \om {\omega}
\def \o{\omega}
\def \s {\sigma}
\def \t {\theta}
\def \z {\zeta}
\def \vp {\varphi}
\def \vt {\vartheta}
\def \ve {\varepsilon}
\def \i {\infty}
\def \ds {\displaystyle}
\def \oo {\overline{\Omega}}
\def \ov {\overline}
\def \bd {\bigtriangledown}
\def \U {\bigcup}
\def \un {\underline}
\def \h {\hspace{.5cm}}
\def \hs {\hspace{2.5cm}}
\def \v {\vspace{.5cm}}
\def \mi {M_{i}}
\def \ra {\longrightarrow}
\def \Ra {\Longrightarrow}
\def \rw {\rightarrow}
\def \bs {\backslash}
\def \rn {{\bf R}^n}
\def \h* {\hspace*{1cm}}
\def\la{\langle}
\def\ra{\rangle}
\def\cal{\mathcal}

\title[Dirac-Harmonic Maps]{Dirac-Harmonic Maps}

\author{Qun Chen, J\"urgen Jost, Jiayu Li, Guofang Wang}

\thanks{The research of QC and JYL was partially supported by  NSFC. QC also thanks the Max
  Planck Institute for Mathematics in the Sciences for support and
  good working conditions during the preparation of this paper. }

\address{School of
Mathematics and Statistics\\ Central China Normal University\\
Wuhan 430079, China } \email{qunchen@mail.ccnu.edu.cn}

\address{Max Planck Institute for Mathematics in the Sciences\\Inselstr. 22-26\\D-04103 Leipzig, Germany
} \email{jjost@mis.mpg.de}

\address{Partner Group of Max Planck Institute for  Mathematics in the Sciences\\
Institute of Mathematics\\
Chinese Academy of Sciences\\
Beijing 100080, P. R. of China} \email{lijia@mail.amss.ac.cn}

\address{Max Planck Institute for Mathematics in the Sciences\\Inselstr. 22-26\\D-04103 Leipzig, Germany
} \email{gwang@mis.mpg.de}

\vskip-15mm
 \vskip-30mm
\today 

\begin{abstract}
We introduce a functional that couples the nonlinear sigma
model with a  spinor field:
   $L=\int_M[|d\phi|^2+(\psi,\D\psi)]$. In two dimensions, it  is
   conformally invariant. The critical points of this 
functional   are called Dirac-harmonic maps. 
  We study some geometric and analytic aspects of such maps, in
  particular a removable singularity theorem.
\end{abstract}
\maketitle

\section{Introduction}
This paper introduces and studies an extension of an established
mathematical subject, namely harmonic maps from Riemann surfaces
into Riemannian manifolds, that is motivated by a model from
quantum field theory, the supersymmetric non-linear sigma model.
What distinguishes harmonic maps from surfaces from those from
higher dimensional domains is the feature of conformal invariance.
On one hand, this is a global aspect, and it implies that when the
domain is the 2-sphere, any such map is automatically conformal;
thus, the second order equations that characterize a harmonic map
reduce to first order equations. As is typical in the geometric
calculus of variations, this leads to important invariants; here,
we obtain minimal surfaces in Riemannian manifolds as well as
quantum cohomology, and the theory of pseudo-holomorphic curves in
symplectic geometry shares the same root. On the other hand, this
leads to a non-compact local invariance group and turns the
regularity and existence problem into a borderline case of the
Palais-Smale condition. In fact, this conformal invariance
connects the local regularity theory with global solutions defined
on the 2-sphere as we know from the seminal paper of
Sacks-Uhlenbeck \cite{SU}. The approach of Sacks-Uhlenbeck depended
on perturbing the functional to ones that satisfy the Palais-Smale
condition and then pass to a limit. This, together with a
sophisticated local analysis, then yielded global existence
results, with possible obstructions coming from the second
homotopy group of the target. (In the absence of those
obstructions, the existence had been shown independently by
Lemaire \cite{Le}.) Subsequently, alternative existence schemes
were developed by Struwe\cite{St} and Chang \cite{Ch} (heat flow
method) and Jost \cite{J1} (local iteration technique). A crucial
point of the work of Sacks-Uhlenbeck was the removability of
isolated singularities through a detailed blow-up analysis. In the
variational context, it is not difficult to reduce the general
regularity question to the one of isolated singularities, and
therefore, the analysis of Sacks-Uhlenbeck was sufficient for the
variational existence scheme. The more general regularity question
for weak solutions of the harmonic map equation was solved later
by H\'elein \cite{He}.\par Now, in the physics literature, the same
model goes under the name non-linear sigma model; when the target
is an $(N-1)$-dimensional sphere, a case of special interest for
quantum field theory, it is called more precisely the non-linear
$O(N)$ sigma model. Here, our map $\phi$ becomes a Bosonic scalar
field  satisfying a non-linear constraint. Now, this model admits
a supersymmetric extension (see \cite{QFT1} for a detailed
exposition\footnote{For a variant of this model, see
  \cite{Sin}}) where $\phi$ is coupled to a
Fermionic field $\psi$. That Fermionic field, and then also the
Bosonic field, is Grassmann valued, a property needed for an
additional symmetry exchanging the two types of fields, called
supersymmetry. One can then try to solve the corresponding
equations by expansions in the degree of the Grassmann
variables.\par Here, we rather let the field $\psi$ also be
real-valued; we then cannot use such expansions anymore, but
rather get a coupled set of field equations for the two fields
$\phi$ and $\psi$. Mathematically, $\phi$ as before is a mapping
from a Riemann surface into a Riemannian manifold, whereas $\psi$
is a spinor field with values in the pull-back of the tangent
bundle of the target under $\phi$. Thus, with respect to
coordinate transformations on the target, $\psi$ transforms as a
tangent vector. The important point is that the Lagrangian is set
up in such a manner as to still be conformally invariant. We thus
obtain an extension of the harmonic map model within the class of
conformally invariant variational problems with a rich internal
structure. Here, we study the geometric setting of the problem and
present a complete variational analysis, that is, we prove the
removability of isolated singularities and the regularity of
solutions arising from a variational scheme. Our analysis is based
on the scheme of Sacks-Uhlenbeck \cite{SU}, but because of the
coupling between the two fields, new difficulties and subtleties
arise that make the analysis considerably harder. A regularity
result in the sense of H\'elein \cite{He} will be presented
elsewhere \cite{CJLW}. Having thus laid the analytical foundation,
the existence question and the construction of geometric
invariants from the solution spaces can be addressed in subsequent
work. We should remark that the analysis in our framework is more
difficult than in the supersymmetric one because we cannot use the
expansions of the physicists. On the other hand, our present model
seems to be very natural from the point of view of Riemannian
geometry, and it thus falls in a similar category as the
Seiberg-Witten or Chern-Simons-Higgs models, and this is our main
motivation for its study. \par\vskip12pt
 Let us now describe the
mathematical framework in more detail: Let $(M, h)$ be an
oriented, compact Riemannian surface and $P_{SO(2)}\to \Sigma$ its
oriented orthonormal frame bundle. A $Spin$-structure is a lift of
the structure group SO(2) to $Spin(2)$, {\it i.e.} there exists a
principal $Spin$-bundle $P_{Spin(2)}\to M$ such that there is a
bundle map
\[ \begin{array}{ccc}
 P_{Spin(2)} & \longrightarrow & P_{SO(2)} \cr
&& \cr \downarrow & & \downarrow \cr && \cr M & \longrightarrow &
M.\cr
\end{array}\]

Let $\Sigma ^+M:=P_{Spin(2)}\times_{\rho}{\mathbb C}$ be a complex
line bundle over $M$ associated to $ P_{Spin(2)}$ and to the
standard representation $\rho:\S^1\to U(1)$. This is the bundle of
positive half-spinors. Its complex conjugate $\Sigma^-M:=
\overline {\Sigma^+M}$ is called the bundle of negative
half-spinors. The spinor bundle is $\Sigma M:=\Sigma ^+M\oplus
\Sigma ^-M.$

There exists a Clifford multiplication
\begin{eqnarray*}
TX\times_{\C} \Sigma^+M &\to &  \Sigma^-M \cr TX\times_{\C}
\Sigma^-M &\to &  \Sigma^+M \end{eqnarray*} denoted by $v\otimes
\psi \to v \cdot \psi$, which satisfies the Clifford relations
\[ v\cdot w \cdot\psi + w\cdot v\cdot \psi =-2h(v,w) \psi,\]
for all $v,w \in TM$ and $\psi \in \Sigma M$.

On the spinor bundle $\Sigma M$ there is a hermitian metric
$ \la \cdot,\cdot \ra $. Let $\n$ be the Levi-Civita connection on $M$
with respect to $h$. There is a  connection (also denoted by $\n$)
on  $\Sigma M$ compatible with the hermitian metric.

Let $\phi$ be a smooth map from $M$ to another Riemannian manifold
$(N,g)$ of dimension $n\ge 2$. Denote $\phi^{-1} TN$ the pull-back
bundle of $TN$ by $\phi$ and consider the twisted bundle $\Sigma
M\otimes \phi^{-1} TN$. On $\Sigma M\otimes \phi^{-1} TN$ there is
a  metric induced from the metrics on $\Sigma M$ and $\phi^{-1}
TN$. Also we have a natural connection $\widetilde{\n}$ on $\Sigma
M\otimes \phi^{-1} TN$ induced from those on $\Sigma M$ and
$\phi^{-1} TN$. In local coordinates, the section $\psi$ of
$\Sigma M\otimes \phi^{-1} TN$ can be expressed by
\[\psi(x)= \sum_{j=1}^n\psi^j (x)\frac{\partial}{\partial y^j}(\phi(x)),\]
where $\psi^i$ is a spinor and $\{\frac{\partial}{\partial y^j}\}$
is the natural local basis. $\widetilde{\n}$ can be expressed by

\[\widetilde{\nabla} \psi= \sum_{i=1}^n\n \psi^i(x) \frac{\partial}{\partial y^j}(\phi(x))
+ \sum_{i,j,k=1}^n\Gamma^i_{jk} \n \phi^j(x)
\psi^k(x)\frac{\partial}{\partial y^i} (\phi(x)).\]

 It is
easy to check that
\[v(\psi_1,\psi_2)=
( \widetilde{\n}_v\psi_1,\psi_2)+
(\psi_1,\widetilde{\n}_v\psi_2),\] for any vector field $v$.

Now we define  the {\it Dirac operator along the map $\phi$} by
\begin{equation}
\label{1.1} \D\psi =\sum_{i}
\partial \hskip -2.2mm \slash
 \psi^i(x)\frac{\partial}{\partial y^i}(\phi(x))
+ \sum_{i,j,k=1}^n\Gamma^i_{jk} \n_{e_\a} \phi^j(x) e_\a\cdot
\psi^k(x)\frac{\partial}{\partial y^i} (\phi(x)),
\end{equation}
 where
$e_1,e_2$ is the local orthonormal basis of $M$ and $
\partial \hskip -2.2mm \slash
:=\sum_{\a=1}^2e_{\a}\cdot \n_{e_\a}$ is the usual Dirac operator.
The Dirac operator $\D$ is formally self-adjoint, i.e.,
\begin{equation}
\label{1.2}
 \int_M(\psi, \D\xi)=\int_M(\D\psi,\xi),
 \end{equation}
for all $\psi, \xi \in \Gamma (\Sigma M\otimes \phi^{-1}TN)$, the
space of smooth section of $\Sigma M\otimes \phi^{-1}TN$. Set
\[{\cal X}:=\{ (\phi,\psi)\,|\, \phi \in C^\infty(M,N) \hbox{ and }
\phi \in \Gamma (\Sigma M\otimes \phi^{-1}TN)\}.\] On $\cal X$, we
consider the following functional
\begin{eqnarray}
\label{1.3}
 L(\phi,\psi)&=&\int_M[|d\phi|^2+(\psi,\D\psi)]\sqrt{h}d^2x \nonumber \\
&=&\int_M [g_{ij}(\phi)h^{\alpha\beta}\frac{\partial
\phi^i}{\partial x_\alpha}\frac{\partial \phi^j}{\partial
x_\beta}+g_{ij}(\phi)(\psi^i,\D\psi^j)]\sqrt{h}d^2x,
\end{eqnarray}
where $h:={\rm det}(h_{\alpha\beta})$. The Euler-Lagrange
equations of $L$ are:
\begin{equation} \label{1.4}
\tau(\phi)={\cal R}(\phi,\psi),
\end{equation}
\begin{equation}
\label{1.5} \D\psi=0,
\end{equation}
where $\tau(\phi)$ is the tension field of the map $\phi$ and
${\cal R}(\phi,\psi) \in \Gamma(\phi^{-1}TM)$ defined by
\begin{equation}
\label{1.5.1} {\cal R}(\phi,\psi)(x) = \frac{1}{2}\sum
R^m_{lij}(\phi(x)) \la \psi^i,\nabla\phi^l\cdot\psi^j \ra
\frac{\partial}{\partial y^m}(\phi(x)).\end{equation} Here
$R^m_{lij}$ are the components of the Riemannian curvature tensor of
$h$. Note that the product $\langle \cdot, \cdot \rangle$
here is the one of $\Sigma M$.
 Solutions $(\phi,\psi)$ to (\ref{1.4}) and (\ref{1.5}) are
called {\it Dirac-harmonic maps}.
\par
As already mentioned, the functional $L$ arises from our study of
the supersymmetric nonlinear sigma model, the difference being
that here $\psi$ is an ordinary (twisted) spinor. Equations
(\ref{1.4}) and (\ref{1.5}) couple the harmonic equation and the
Dirac equation in a rather natural way.

In this paper, our main aim is to establish some geometric and
analytic aspects of solutions to equations (\ref{1.4}) and
(\ref{1.5}). We first establish some basic properties of the
Dirac-harmonic maps and give some examples of nontrivial
solutions. We also derive some geometric properties of
 the Dirac-harmonic maps,  the conformal invariance and
 the existence of a generalized Hopf differential.
 Then we prove the removablity of singularities for the solutions.
\par This paper is organized as follows: in Section 2, we deduce
the Euler-Lagrange equations, and construct nontrivial solutions;
in Section 3, we define the energy-momentum tensor of the action
$L$ and construct a holomorphic differential (Proposition 3.2)
which plays a role in proving the removable singularity theorem;
we also establish the basic Weitzenb\"{o}ck formula (Proposition
3.4) for spinor fields satisfying (\ref{1.5}); in Section 4 we
 prove the main result about removable singularities (Theorem 4.7).

\section{Dirac-harmonic maps}
\addtocounter{equation}{-6}

In this section, we establish some basic facts for the functional
$L$ and equations (\ref{1.4})--(\ref{1.5}). 
\vskip12pt
\par
\noindent {\bf Proposition 2.1.} {\it The Euler-Lagrange equations
for $ L$ are
\begin{eqnarray}
 \label{2.1.2} \tau(\phi)&=&{\cal
R}(\phi,\psi) \\
\label{2.1.1} \D \psi &=&0,
\end{eqnarray}
where $\tau(\phi)$ is the tension field of the map $\phi$ and
${\cal R}$ is defined by (\ref{1.5.1}).}

\vskip12pt

 \pr Equation (\ref{2.1.1}) is easy to derive. Consider
 a family of $\psi_t$ with $d \psi_t / dt=\eta$ at $t=0$
 and fix $\phi$.
 Since $\D$ is
formally self-adjoint, we have
\begin{eqnarray*}
\frac{dL}{dt}|_{t=0} &=&
=\int_M \la \eta,\D\psi \ra +\la \psi,\D\eta\ra\\
&=&2\int_M \la \eta,\D\psi \ra.
\end{eqnarray*}
Hence, we get (\ref{2.1.1}).
\par
Next, we consider a variation $\{ \phi_t\}$ of $\phi$ such that
$d\phi_t/dt=\xi$ at $t=0$ and fix $\psi$. We choose $\{e_\alpha\}$ as 
a local orthonormal basis on $M$ such that $[e_\alpha,
\partial_t]=0$, $\nabla_{e_\alpha}e_\beta=0$
 at a considered point.
\begin{equation}
\label{2.1.2.0}
\frac{dL(\phi_t)}{dt}|_{t=0}=\int_M\frac{\partial}{\partial
t}|d\phi_t|^2|_{t=0}+\int_M\frac{\partial}{\partial t} \la
\psi,\D\psi \ra |_{t=0}:=I+II.
\end{equation}
It is easy to check that
\begin{equation}
\label{2.1.2.1} I=-2\int_M\tau^i(\phi)g_{im}\xi^m.
\end{equation}
See the proof for instance in \cite{Xin}. Now we compute II. First
we compute the variation of $\D \psi$. We have
\begin{eqnarray*}
\frac{d}{dt}\D\psi&=&e_\alpha\cdot\nabla_{\partial_t}\nabla_{e_\alpha}\psi\\
&=&
e_\alpha\cdot\nabla_{e_\alpha}\psi^i\otimes\nabla_{\partial_t}\partial_{y_i}+e_\alpha\cdot
\psi^i\otimes\nabla_{\partial_t}\nabla_{e_\alpha}\partial_{y_i}\\
&=&e_\alpha\cdot\nabla_{e_\alpha}\psi^i\otimes\nabla_{\partial_t}\partial_{y_i}+e_\alpha\cdot
\psi^i\otimes[\nabla_{e_\alpha}\nabla_{\partial_t}\partial_{y_i}+R(\partial_t,e_\alpha)\partial_{y_i}]\\
&=&
e_\alpha\cdot\nabla_{e_\alpha}(\psi^i\otimes\nabla_{\partial_t}\partial_{y_i})+e_\alpha\cdot
\psi^i\otimes R^N(d\phi(\partial_t),
d\phi(e_\alpha))\partial_{y_i}.
\end{eqnarray*}
Hence, we have
\begin{eqnarray*}
II&=&\int_M \la \xi,\D\psi \ra +
\int_M \la \psi,\frac{d}{dt}\D\psi\ra|_{t=0}\\
&=& \int_M \la
\psi,\D(\psi^i\otimes\nabla_{\partial_t}\partial_{y_i}))|_{t=0}+ \la\psi,e_\alpha\cdot
\psi^i\otimes R^N(d\phi(\partial_t),
d\phi(e_\alpha))\partial_{y_i} \ra |_{t=0}\\
&=&\int_M \la
\D\psi,\psi^i\otimes\nabla_{\partial_t}\partial_{y_i}\ra |_{t=0}+
\la \psi,e_\alpha\cdot \psi^i\otimes R^N(d\phi(\partial_t),
d\phi(e_\alpha))\partial_{y_i} \ra
|_{t=0}\\
&=&\int_M \la \psi,e_\alpha\cdot \psi^i\otimes
R^N(d\phi(\partial_t),
d\phi(e_\alpha))\partial_{y_i} \ra |_{t=0}\\
&=&\int_M \la \psi,e_\alpha\cdot \psi^i\otimes
R^N(\xi^m\partial_{y_m},\phi^l_\alpha
\partial_{y_l})\partial_{y_i} \ra \\
&=&\int_M \la \psi,e_\alpha\cdot \psi^i\otimes \xi^m\phi^l_\alpha
R^j_{iml}\partial_{y_j}\ra \\
&=&\int_M \la \psi^i,\nabla\phi^l\cdot\psi^j \ra R_{mlij}\xi^m,
\end{eqnarray*}
where we have used (\ref{2.1.1}).
Consequently, we have
$$\frac{dL(\phi_t)}{dt}|_{t=0}=\int_M[-2g_{mi}\tau^i(\phi)+R_{mlij}
(\psi^i,\nabla\phi^l\cdot\psi^j)]\xi^m,$$ and hence (\ref{2.1.2}).
 \qed
\vskip12pt

It is obvious that there are two types of trivial solutions. One
is $(\phi, 0)$, where $\phi$ is a harmonic map, and another is
$(y, \psi)$, where $y$ is a point in $N$ viewed as a constant map
from $M\to N$, and $\psi$ is a harmonic spinor, i.e., $\partial
\hskip -2.2mm \slash\psi=0$.

Let us give a construction of non-trivial solutions. Let $M=\S^2$
and $N=\S^2$. Let $\Sigma \S^2$ be the spinor bundle over $\S^2$
with respect to the unique spin structure. For any map
$\phi:\S^2\to \S^2$ and a spinor $\Psi \in \Gamma(\Sigma \S^2)$,
we define a spinor field $\psi$ along the map $\phi$ by
\begin{equation}
\label{2.1.3} \psi_{\phi,\Psi} = e_\alpha\cdot \Psi
\phi_*(e_\alpha),\end{equation} where $e_\alpha$ ($\alpha=1,2$) is
a local basis of $T\S^2$. It is clear that $\psi_{\phi,\Psi}$ is
well-defined. A spinor $\Psi \in \Gamma (\Sigma \S^2)$ is called a
 twistor spinor if
\[ \n _v \Psi+\frac 1 2 v\cdot \partial \hskip -2.2mm \slash \Psi =0,\]
for any vector field $v\in \Gamma (T\S^2)$.

 Now we have
\vskip12pt
 \noindent
 {\bf Proposition 2.2.} {\it
Let $\psi_{\phi,\Psi}$ be defined by (\ref{2.1.3}) from a
nonconstant map $\phi: \S^2\to \S^2$ and a spinor $\Psi$. Then
$(\phi, \psi_{\phi,\Psi})$ is a Dirac harmonic map if and only if
$\phi$ is a (possibly branched) conformal map and $\Psi$ is a twistor spinor.
} \vskip12pt

\pr We first assume that  $(\phi,\psi_{\phi,\Psi})$ is a
Dirac-harmonic map. For the spinor field $\psi$ constructed from
 (\ref{2.1.3}), we always have
\begin{eqnarray*}
\la \psi^k,\nabla \phi^j\cdot\psi^l \ra  &=& \la \nabla
\phi^k\cdot\Psi,\nabla
\phi^j\cdot\nabla\phi^l\cdot\Psi \ra \\
&=&\phi^k_\alpha\phi^j_\beta\phi^l_\gamma\la
e_\alpha\cdot\Psi,e_\beta\cdot
 e_\gamma\cdot\Psi \ra \\
&=&0 ,
\end{eqnarray*}
since $ e_\beta\cdot e_\gamma=-e_\gamma\cdot e_\beta $ for $
\beta\not=\gamma $ and $\phi_\beta\phi_\gamma$ is symmetric.
Consequently,
$$\frac 1 2 R^i_{jkl}(\psi^k,\nabla\phi^j\cdot\psi^l)\equiv 0.$$
 Hence, if $(\phi, \psi_{\phi,\Psi})$ is a Dirac harmonic
map, then $\phi$ is harmonic from equation (\ref{1.4}). Choosing
suitable coordinates such that $\n _{e_\alpha}e_\beta=0$ and $x\in
M$, by (\ref{1.5}) at $x$ we have
\begin{eqnarray*}
0&=& \D \psi_{\phi,\Psi}\\
&=& e_\beta\cdot \widetilde{\n}_{e_\beta}(e_\alpha\cdot \Psi \phi_*(e_\alpha))\\
&=& e_\beta\cdot e_\alpha\cdot\{\n _{e_\beta}\Psi
\phi_*(e_\alpha)+
\Psi \n_{e_\beta}\phi_*(e_\alpha)\}\\
&=&-\{\n _{e_\alpha}\Psi \phi_*(e_\alpha)+\Psi \tau(\phi)\}\\
&&+e_1\cdot e_2\cdot\{\n _{e_1}\Psi \phi_*(e_2)-\n _{e_2}\Psi
\phi_*(e_1)
+\Psi [\n_{e_1}\phi_*(e_2)-\n_{e_1}\phi_*(e_2) ]\}\\
&=&-\n _{e_\alpha}\Psi \phi_*(e_\alpha) -e_1\cdot e_2\cdot\{\n
_{e_2}\Psi \phi_*(e_1)-\n _{e_1}\Psi \phi_*(e_2)\}
\end{eqnarray*}
Since $\phi$  is conformal (and non-constant), the above equation is
equivalent to 
\begin{equation}
\label{eq2.1} e_1\cdot \n _{e_1}\Psi=e_2\cdot \n _{e_2}\Psi,
\end{equation}
 which is equivalent to the condition
 that $\Psi$ is a twistor spinor.
 \par
Conversely, from the above proof, it is easy to see that if $\phi$
is a conformal map and $\Psi$ is a twistor spinor, then $(
\phi,\psi_{\phi,\Psi}) $ is a Dirac-harmonic map. \qed \vskip12pt
\noindent {\bf Remark 2.3.} The twistor spinors form a linear
space of complex dimension 2, whereas the conformal maps form a
nonlinear space whose dimension depends on their degree; for
degree 1, we have dimension 3.

\section{Geometric aspects of Dirac-harmonic maps}
\addtocounter{equation}{-6}

In this section, we will study geometric properties of
Dirac-harmonic maps. First
 \vskip12pt \noindent {\bf Lemma
3.1.} {\it
The functional $L(\phi,\psi)$ is conformally invariant. Namely,
for any conformal diffeomorphism $f:M\to M$, set \[\tilde
\phi=\phi \circ f \quad \hbox{ and }\quad \tilde \psi = \l^{-1/2}
\psi\circ f.\] Then $L(\phi,\psi)=L(\tilde \phi,\tilde\psi).$ Here
$\l$ is the conformal factor of the conformal map $f$.
} \vskip12pt \pr The proof is standard. In fact, the following
terms are invariant under a conformal transformation
\[\int |\n \phi|^2 dvol(g), \int \la \psi, \D\psi  \ra dvol(g)
\hbox{ and }\int |\psi|^4 dvol(g).\] Here we only check the
conformality of $\int \la \psi, \D\psi \ra dvol(g)$. Let $\tilde
g=f^*g$ and $\widetilde \partial \hskip-2.2mm \slash$ the usual
Dirac operator with respect to the new metric $\tilde g$. By the
conformality of $f$, we have $\tilde g=(\l^*)g$. We identify the
new and old spin bundles as in \cite{Hi}. Recall that the relation
between two Dirac operators $\partial \hskip-2.2mm \slash$ and
$\widetilde \partial \hskip-2.2mm \slash$ is (cf. \cite{Hi})
\begin{equation}
\label{2.7} \widetilde{\partial \hskip-2.2mm
\slash}\widetilde{\psi}=\lambda^{-\frac{3}{2}}\partial
\hskip-2.2mm \slash \psi,
\end{equation}
from which one can deduce that
\begin{equation}
\label{2.8}
\widetilde{\D}\widetilde{\psi}=\lambda^{-\frac{3}{2}}\D \psi,
\end{equation}
and hence the conformal invariance of $\int \la \psi, \D\psi\ra
dvol(g)$. \qed
 \vskip12pt
For a two dimensional harmonic map, there is an important
holomorphic quadratic differential, the Hopf differential. For a
Dirac-harmonic map, we also have an analogue.

Let $(\phi,\psi)$ be a Dirac-harmonic map.  On a small domain
$\Omega$ of $M$, choose $z=x+{ i}y$ a local isothermal
parameter $z=x+iy$ with $g=ds^2=\rho |dz|^2$. Define
\begin{equation} \label{3.8} T(z){
dz}^2=\{(|\phi_x|^2-|\phi_y|^2-2i\la\phi_x,\phi_y\ra)+
(\la\psi,\partial_x\cdot \widetilde{\nabla}_{\partial_x}\psi \ra
-{ i} \la \psi,\partial_x\cdot
\widetilde{\nabla}_{\partial_y}\psi\ra )\}{ dz}^2.
\end{equation}
Here $\partial_x=\frac{\partial}{\partial x}$ and
$\partial_y=\frac{\partial}{\partial y}$. \vskip12pt \noindent
 {\bf Proposition 3.2.}
 {\it The quadratic differential $Tdz^2$ is holomorphic.}
\par
\vskip12pt One can prove this proposition by a direct computation,
or as a consequence of a conservation law as follows.

Define a two-tensor by
\begin{equation}
\label{3.1} T_{\alpha\beta}:=(2 \la \phi_\alpha,\phi_\beta \ra-
\delta_{\alpha\beta} \la \phi_\gamma,\phi_\gamma\ra)+
\la\psi,e_\alpha\cdot\wb\psi\ra,
\end{equation}
where $\phi_\alpha:=\phi_*(e_\alpha).$ Here,  $\{e_\alpha\}$ is a
local orthonormal basis on $M$ and $\{\theta^\alpha\}$ a coframe
dual to $\{e_\alpha\}$. The tensor
$T_{\alpha\beta}\theta^\alpha\otimes\theta^\beta$ is called the 
energy-momentum tensor. Using the equation $\D\psi=0$ one can
verify that $T$ is symmetric:
\[T_{\alpha\beta}=T_{\beta\alpha}.\]
 \noindent
 {\bf Proposition 3.3.}
 {\it Let $(\phi,\psi)$ be a smooth solution of (\ref{1.4}) (\ref{1.5})
on $M$, then the energy-momentum tensor is conserved, i.e.,
\begin{equation}\label{3.2}
\sum\limits_\alpha\nabla_{e_\alpha}T_{\alpha\beta}=0.\end{equation}}
\pr
\begin{eqnarray*}
\nabla_{e_\alpha}T_{\alpha\beta}&=&\nabla_{e_\alpha} (2
\la\phi_\alpha, \phi_\beta \ra-\delta_{\alpha\beta} \la
\phi_\gamma,\phi_\gamma \ra )+ \nabla_{e_\alpha}
                                  \la \psi,e_\alpha\cdot\wb\psi \ra\\
                                  &:=& I+II.
\end{eqnarray*}
We choose a local orthonormal basis $\{e_\alpha\}$ on $M$ such
that $\nabla_{e_\alpha}e_\beta=0$ at a considered point. We
compute
\begin{eqnarray*}
I &=& 2 \la
\nabla_{e_\alpha}\phi_*(e_\alpha),\phi_*(e_\beta)\ra+2\la
\phi_*(e_\alpha),\nabla_{e_\alpha}\phi_*(e_\beta) \ra \\
      &&-2\delta_{\alpha\beta}\la\phi_*(e_\gamma),\nabla_{e_\alpha}\phi_*(e_\gamma)\ra\\
      &=&2\la\tau(\phi),\phi_\beta\ra+2\la\phi_\alpha,\nabla_{e_\beta}\phi_*(e_\alpha)\ra-
      2\la\phi_\gamma,\nabla_{e_\beta}\phi_*(e_\gamma)\ra\\
      &=&2\la\tau(\phi),\phi_\beta \ra ,
\end{eqnarray*}
and
\begin{eqnarray*}
II&=& \la \psi_\alpha,e_\alpha\cdot\psi_\beta \ra + \la
\psi,e_\alpha\cdot\wa\wb\psi \ra
\\
&=&- \la e_\alpha\cdot\psi_\alpha,\psi_\beta \ra + \la \psi,\D\psi_\beta\ra \\
&=& \la \psi,\D\psi_\beta \ra.
\end{eqnarray*}
Therefore, we have
\begin{equation}
\label{3.4} \nabla_{e_\alpha}T_{\alpha\beta}=2 \la
\tau(\phi),\phi_\beta \ra+ \la\psi,\D\psi_\beta\ra .
\end{equation}
Now
\begin{eqnarray}
\label{3.5} 2 \la \tau(\phi),\phi_\beta \ra &=&2\la
\frac{1}{2}R^m_{lij} \la\psi^i,\nabla\phi^l\cdot\psi^j\ra
\partial_{y^m},\phi^p_\beta\partial_{y^p} \ra \nonumber
\\
&=& g_{mp}R^m_{lij} \la \psi^i,\nabla\phi^l\cdot\psi^j
\ra\phi^p_\beta \nonumber
\\
&=&R_{mlij} \la \psi^i,\nabla\phi^l\cdot\psi^j\ra\phi^m_\beta.
\end{eqnarray}
We compute $\D\psi_\beta=e_\alpha\cdot\wa\wb\psi.$ By a direct
computation, we have:
$$\wa\wb\psi-\wb\wa\psi=R^{\Sigma M}(e_\alpha,e_\beta)
\psi^i\otimes\partial_{y^i}+R^m_{lij}\phi^i_\alpha
\phi^j_\beta\psi^l\otimes\partial_{y^m},$$ where $R^{\Sigma M}$ is
the curvature operator of the connection $\n$ on the spinor bundle
$\Sigma M$. For this  curvature operator, we have (cf. \cite{J})
\begin{equation}
\label{Ric} e_\alpha\cdot R^{\Sigma
M}(e_\alpha,X)\psi^i=\frac{1}{2}Ric(X)\cdot\psi^i, \quad \forall
X\in \Gamma(M).\end{equation} It follows that
\begin{eqnarray*}
\la \psi,e_\alpha\cdot R^{\Sigma M}(e_\alpha,e_\beta)
\psi^i\otimes\partial_{y^i} \ra &=& \la
\psi^j\otimes\partial_{y^j},e_\alpha\cdot R^{\Sigma
M}(e_\alpha,e_\beta)
\psi^i\otimes\partial_{y^i} \ra\\
&=&g_{ij} \la \psi^j,e_\alpha\cdot R^{\Sigma M}(e_\alpha,e_\beta)\psi^i\ra \\
&=&\frac{1}{2}g_{ij} \la \psi^j,Ric(e_\beta)\cdot\psi^i\ra\\
&=& 0.
\end{eqnarray*}
Therefore
\begin{equation}\begin{array}{rcl}
\label{3.6} \la\psi,\D\psi_\beta\ra&=&
\la\psi,e_\alpha\cdot\wa\wb\psi \ra \\
&=&
\la \psi,\wb(e_\alpha\cdot\wa\psi) \ra +R^m_{lij}\phi^j_\beta \la
\nabla\phi^i\cdot\psi^l\otimes\partial_{y^m},\psi^p\otimes\partial_{y^p}
\ra\nonumber
\\
&=&R^m_{lij}\phi^j_\beta \la \psi^p,\nabla\phi^i\cdot
\psi^l\ra g_{mp} \\
&=&-R_{mlij} \la \psi^i,\nabla\phi^l\cdot \psi^j\ra \phi^m_\beta.
\end{array}\end{equation}
  From (\ref{3.4}),(\ref{3.5}) and (\ref{3.6}) we conclude that
$T_{\alpha\beta}$ is conserved. \qed
\par
\vskip12pt \noindent
{\it Proof of Proposition 3.2.} The proof
follows directly from Proposition 3.3. \qed
\medskip
\par
We now consider  (\ref{1.5}), i.e.,
$$\D\psi=(
\partial \hskip -2.2mm
\slash
\psi^i+\Gamma^i_{jk}(\phi)\partial_\alpha\phi^je_\alpha\cdot\psi^k)
\otimes\partial_{y^i}=0,$$ for  $\phi: M^m \to N^n,$ $\psi\in
\Gamma(\Sigma M\otimes\phi^{-1}TN)$.
Here $M^n$ is an $n$ dimensional spin manifold with $n\ge 2$.
Clearly, one can also discuss (1.4)-(1.5) for higher dimensional spin
manifolds.
 We  call such a $\psi$ a
harmonic spinor field along the map $\phi$. Note that here we do
not assume $(\phi,\psi)$ to be a solution of (\ref{1.4}).

For a harmonic spinor field along a map $\phi$, we have the
following Weitzenb\"{o}ck formula. \vskip12pt \noindent {\bf
Proposition 3.4.}
 {\it
 Let $M^m$ and $N^n$ be Riemannian manifolds, $\phi : M\to N$,
 and
 $\psi=\psi^i\otimes\frac{\partial}{\partial y^i}(\phi)\in \Gamma (\Sigma M\otimes
\phi^{-1}TN)$. Then
\begin{eqnarray}
 \label{3.1.3}
\D\hskip1mm^2\psi &= &-\wa\wa\psi +\frac{1}{4}
R\psi+\frac{1}{2}R^i_{kpj}\phi^p_\alpha\phi^j_\beta(e_\alpha\cdot
e_\beta\cdot\psi^k)\otimes\frac{\partial}{\partial y^i}.
\end{eqnarray}
If $\phi$ and $\psi$ satisfy (\ref{1.5}),
then
\begin{equation}
\label{4.5} \frac{1}{2}\Delta |\psi
|^2=|\widetilde{\nabla}\psi|^2+\frac{1}{4}R|\psi|^2-\frac{1}{2}R_{ijkl}(\nabla
\phi^k\cdot \psi^i,\nabla\phi^l\cdot \psi^j),
\end{equation}
where $R$ is the scalar curvature of  $M$, $\widetilde{\nabla}$
denotes the connection on $\Sigma M\otimes \phi^{-1}TN$ and
$|\psi|^2:=g_{ij}(\phi)(\psi^i,\psi^j)$. } \vskip 12pt \noindent
{\it Proof.} One can apply a general  Weitzenb\"ock formula (see for
example \cite{LM} or \cite{J}) to
prove the proposition. For the convenience of the reader, we
present a proof here.
 Choose an orthonormal basis
$\{e_\alpha|\alpha=1,2,\cdots,m\}$
 on $M$ such that $\nabla_{e_\beta}e_\alpha=0$ at a considered
 point.
 Noting that $ e_\alpha^2=-1$, $e_\alpha e_\beta +e_\beta
 e_\alpha =0$ and $\alpha\not=\beta,$
 we have
 \begin{eqnarray}
 \label{3.1.1}
\D\hskip1mm^2\psi &=&  e_\beta\cdot e_\alpha\cdot \wb \wa \psi
\\&= &-\wa\wa\psi +\sum\limits_{\alpha<\beta}
e_\alpha\cdot
e_\beta\cdot \widetilde{R}(e_\alpha,e_\beta)\psi  \nonumber \\
&=& -\wa\wa\psi+\frac{1}{2}\sum\limits_{\alpha,\beta}
e_\alpha\cdot e_\beta\cdot \widetilde{R}(e_\alpha,e_\beta)\psi
\end{eqnarray}
where $\widetilde{R}(\cdot,\cdot)$ is the curvature operator on
$\Sigma M\otimes\phi^{-1}TN$, namely,
$$\widetilde{R}(e_\alpha,e_\beta)\psi=\wa\wb\psi-\wb\wa\psi-\widetilde{\nabla}_{[e_\alpha,e_\beta]}\psi.$$
Since
$$\wa\wb\psi=(\nabla_{e_\alpha}\nabla_{e_\beta}\psi^i+\Gamma^i_{jk,p}\phi^p_\alpha\phi^j_\beta\psi^k)\otimes\frac{\partial}{\partial
y^i} $$ and
$$\wb\wa\psi=(\nabla_{e_\beta}\nabla_{e_\alpha}\psi^i+\Gamma^i_{jk,p}\phi^p_\beta\phi^j_\alpha\psi^k)\otimes\frac{\partial}{\partial
y^i}, $$ we have
\begin{eqnarray*}
\widetilde{R}(e_\alpha,e_\beta)\psi &=& R^{\Sigma
M}(e_\alpha,e_\beta)\psi^i\otimes\frac{\partial}{\partial
y^i}+(\Gamma^i_{jk,p}-\Gamma^i_{kp,j})\phi^p_\alpha\phi^j_\beta\psi^k\otimes\frac{\partial}{\partial
y^i} \\
&=& R^{\Sigma
M}(e_\alpha,e_\beta)\psi+R^i_{kpj}\phi^p_\alpha\phi^j_\beta\psi^k\otimes\frac{\partial}{\partial
y^i}.
\end{eqnarray*}
Putting this formula into (\ref{3.1.1}) we have
\begin{eqnarray}
 \label{3.1.2}
\D\hskip1mm^2\psi &= &-\wa\wa\psi +\frac{1}{2}e_\alpha\cdot
e_\beta\cdot
R^{\Sigma M}(e_\alpha,e_\beta)\psi  \nonumber \\
&& +\frac{1}{2}R^i_{kpj}\phi^p_\alpha\phi^j_\beta(e_\alpha\cdot
e_\beta\cdot\psi^k)\otimes\frac{\partial}{\partial y^i},
\end{eqnarray}
where $R^{\Sigma M}(\cdot,\cdot)$ is the curvature operator on
$\Sigma M$. It is known (cf. \cite{J}) that
$$e_\alpha\cdot e_\beta\cdot
R^{\Sigma M}(e_\alpha,e_\beta)\psi^i=\frac{1}{2}R\psi^i.$$ Thus,
we obtain
\begin{eqnarray}
\D\hskip1mm^2\psi &= &-\wa\wa\psi +\frac{1}{4}
R\psi+\frac{1}{2}R^i_{kpj}\phi^p_\alpha\phi^j_\beta(e_\alpha\cdot
e_\beta\cdot\psi^k)\otimes\frac{\partial}{\partial y^i},
\end{eqnarray}
from which
$$ \la \D\hskip1mm^2\psi,\psi \ra =- \la \wa\wa\psi,\psi\ra +\frac{1}{4}
R|\psi|^2+\frac{1}{2}R_{ikpj}\phi^p_\alpha\phi^j_\beta(e_\alpha\cdot
e_\beta\cdot\psi^k,\psi^i).$$ Therefore,
\begin{eqnarray*}
\frac{1}{2}\Delta |\psi|^2 &=&
 \la \wa\wa\psi,\psi \ra +\la \wa\psi,\wa\psi\ra  \\
&=& |\widetilde{\nabla}\psi|^2- \la \D\hskip1mm^2\psi,\psi\ra
+\frac{1}{4}
R|\psi|^2-\frac{1}{2}R_{ijkl}(e_\alpha\cdot\psi^i,e_\beta\cdot\psi^j)\phi^k_\alpha\phi^l_\beta.
\end{eqnarray*}
 \hfill $\Box$
\par

\vskip12pt
 Let $(N', g')$ be another Riemannian manifold and $f:
N\to N'$ a smooth map. For any $(\phi, \psi)\in {\cal X}$ we set
\[\phi ' =f\circ\phi \quad \hbox{ and }  \quad  \psi '= f_*\psi.\]
It is clear that $ \psi '$ is a spinor along the map $ \phi '$.
Let $A$ be the second fundamental form of $f$, i.e., $A(X,Y)=(
\n_X df)(Y)$ for any $X,Y\in \Gamma(TN)$. The tension fields of
$\phi$ and $\phi '$ have the following relation
\begin{equation}
\label{3.1.5} \tau ( \phi ')=\sum_{\a=1}^2A(d\phi(e_\a),
d\phi(e_\a))+ df(\tau(\phi)).
\end{equation}
It is also easy to check that the Dirac operators $\D$ and
$\D\hskip1mm '$ corresponding to $\phi$ and $\phi '$ respectively
are related by the following formula
\begin{equation}
\label{3.1.6} \D \hskip1mm ' \psi '= f_*(\D\psi)+A(d\phi(e_\a),
e_\a\cdot \psi ).
\end{equation}
Here in local coordinates
\[A(d\phi(e_\a),e_\alpha\cdot \psi )= \phi^i_\alpha
 A(\frac{\partial}{\partial y^i},
\frac{\partial}{\partial y^j})e_\alpha \cdot \psi^j,\] where
$\{\frac{\partial}{\partial y^i}\}$ is a local basis of $N$.
Furthermore, if $f: N\to N'$ is an isometric immersion, then
$A(\cdot,\cdot)$ is the second fundamental form of the submanifold
$N$ in $N'$, and
\[\n'_X\xi=-P(\xi;X)+\n_X^{\perp}\xi, \qquad \quad \n'_XY=\n_XY+A(X,Y)  \]
$\forall X,Y\in \Gamma(TN),$ $\xi\in\Gamma (T^{\perp}N)$, where
$P(\xi; \cdot)$ denotes the shape operator. Note that $ \la
P(\xi;X), Y\ra=\la A(X,Y), \xi \ra.$ We can rewrite equations
(\ref{1.4})-(\ref{1.5}) in terms of $A$ and the geometric data of
the ambient space $N'$. By the equation of Gauss, one has
\begin{equation}
\label{3.1.7} \cal
R(\phi,\psi)=P(A(d\phi(e_\alpha),e_\alpha\cdot\psi);\psi)+\frac 1
2 R'(e_\alpha\cdot\psi,\psi)d\phi(e_\alpha).
\end{equation}
Therefore, by using (\ref{3.1.5}) and (\ref{3.1.6}), and
identifying $\phi$ with $ \phi '$ and $\psi$ with $ \psi '$, we
can rewrite (\ref{1.4}) and (\ref{1.5})  as follows:
\begin{equation}
\label{3.1.9} \tau ( \phi )=A(d\phi(e_\a), d\phi(e_\a))+ \frac 1 2
R'(e_\alpha\cdot\psi,\psi)d\phi(e_\alpha)+P(A(d\phi(e_\alpha),e_\alpha\cdot\psi);\psi),
\end{equation}
\begin{equation}
\label{3.1.8} \D \hskip1mm ' \psi = A(d\phi(e_\a), e_\a\cdot \psi
).
\end{equation}
In particular, if $N'=\R^K$, then these become
\begin{eqnarray}
\label{3.1.11} -\Delta \phi &=& A(d\phi,
d\phi)+P(A(d\phi(e_\alpha),e_\alpha\cdot\psi);\psi) \\
\label{3.1.10}
\partial \hskip -2.2mm \slash
 \psi &=& A(d\phi(e_\a), e_\a\cdot \psi ).
\end{eqnarray}

\vskip24pt

\section{ Analytic aspects: removable singularities}
\addtocounter{equation}{-22} Embed $(N,h)$ into some $\R^K$
isometrically and denote by $A(\cdot,\cdot)$ the second
fundamental form as in the previous section. Any map $\phi$ from
$(M,g)$ to $(N,h)$ can be seen as a map $\phi$ from $(M,g)$ to
$\R^K$ with $\phi(x)\in N$. And any spinor field $\psi$ along the
map $\phi$ can be seen as a $K$-tuple of (usual) spinors
$(\psi^1,\psi^2,\cdots,\psi^K)$  satisfying the condition that for
any normal vector $\nu=\sum_{i=1}^K \nu_i E_i$ of $N$ at
$\phi(x)$, we have
\[\sum_i\nu_i\psi^i=0,\]
where $\{E_i, i=1,2,\cdots, K\}$ is the standard basis of $\R^K$.
In this section, we always view $(\phi,\psi)$ in this way. By the
discussion in the previous section, such a pair $(\phi,\psi)$ is a
Dirac-harmonic map if and only if $(\phi,\psi)$ satisfies
(\ref{3.1.11}) and (\ref{3.1.10}).
\par In this
section, we will prove the removable singularity theorem for
Dirac-harmonic maps with ``finite energy" defined below.
 \vskip12pt \noindent
{\bf Definition 4.1.}  Let $U$ be a domain on $M$. The energy of
$(\phi,\psi)$ on $U$ is:
\begin{equation}
\label{6.20} E(\phi,\psi,U):=\int_U (|d\phi|^2+|\psi|^4).
\end{equation}

 Note that the energy is conformally invariant.
It is crucial for our results.

 Before we consider the analytic aspects of Dirac-harmonic maps, let us note
that on a surface the (usual) Dirac operator $\partial \hskip
-2.2mm \slash$ can be seen as the  Cauchy-Riemann operator.
Consider $\R^2$ with the Euclidean metric $dx^2+dy^2$. Let
$e_1=\frac{\partial}{\partial x}$ and
$e_2=\frac{\partial}{\partial y}$ be the standard orthonormal
frame. A spinor field is simply a map $\Psi:\R^2\to
\Delta_2=\C^2$, and $e_1$ and $e_2$ acting on  spinor fields can
be identified by multiplication with matrices
\[e_1=\left(\begin{matrix}0& 1\\ -1&0 \end{matrix}\right),
\quad e_2=\left(\begin{matrix}0& i\\ i&0 \end{matrix}\right).\] If
$\Psi:=\left(\begin{matrix} \ds f
\\ \ds g\end{matrix}\right)
:\R^2\to \C^2$ is a spinor field, then the Dirac operator is
\[\partial \hskip
-2.2mm \slash\Psi=\ds \left(\begin{matrix}0& 1\\ -1&0
\end{matrix}\right) \left(\begin{matrix} \ds \frac{\partial
f}{\partial x}\\ \ds \frac{\partial g}{\partial x}
\end{matrix}\right)+
\left(\begin{matrix}0& i\\ i&0 \end{matrix}\right)
\left(\begin{matrix} \ds \frac{\partial f}{\partial y} \\
\ds\frac{\partial g}{\partial y}
\end{matrix}\right)=
2\left(\begin{matrix} \ds \frac{\partial g}{\partial \bar z}
\\ -\ds\frac{\partial f}{\partial z}\end{matrix}\right),\]
where
\[\frac{\partial}{\partial z}=\frac 12 \left(\frac{\partial }{\partial x}
- i\frac{\partial }{\partial y}\right), \quad
\frac{\partial}{\partial \bar z}=\frac 12 \left(\frac{\partial
}{\partial x} + i\frac{\partial }{\partial y}\right).\] Therefore,
the elliptic estimates developed for (anti-) holomorphic functions
can be used to study the Dirac equation.
\par
\medskip
 \noindent{\bf Proposition 4.2} {\it Let
 $(N^n,h)$ be a compact Riemannian manifold and $(M^2, g)$ a surface
  with a fixed
 spin structure. Then there is a small constant $\varepsilon >0$
  such that if $(\phi,\psi)$ is a smooth
   Dirac-harmonic map satisfying
 \begin{equation}
 \label{4.1}
 \int_M(|d\phi|^2+|\psi|^4)  < \varepsilon,
 \end{equation}
then $\phi$ is constant and consequently $\psi$ is a usual
harmonic spinor.} \vskip12pt

\pr In view of (\ref{3.1.11}), we have
$$|\Delta\phi|\leq C(|d\phi|^2+|d\phi||\psi|^2),$$
where $C>0$ is a constant depending only on $N$.  Hence we have
\begin{equation}
\label{4.2} \begin{array}{rcl} \ds
\|\Delta\phi\|_{L^{\frac{4}{3}}} &\leq&\ds
C(\|\n \phi\|^2_{L^{\frac{4}{3}}}+\|\psi\|^2_{L^4}\|\n \phi\|_{L^4})\\
&\le&  \ds
C(\|\n \phi\|_{L^{2}}+\|\psi\|^2_{L^4})\| \n \phi\|_{L^{1,\frac 4 3}}\\
&\le & \varepsilon C \|\n \phi\|_{L^{1,\frac 4 3}}
\end{array}\end{equation}
Therefore, if $\varepsilon$ is small enough, we can show that  $\phi\equiv
const$, and hence $\psi$ is a harmonic spinor.
 \hfill $\Box$
\par
Now we consider the local behavior of Dirac-harmonic maps. Since
they
 are conformally invariant,
in the sequel we may assume $M$ to be the unit disk $D$ with
trivial spin structure.
\par
\vskip12pt
 \noindent {\bf Theorem 4.3.} {\it There is a small constant $\varepsilon
>0$ such that if $(\phi,\psi)$ is 
a Dirac harmonic map satisfying
\begin{equation}
 \label{5.2}
 \int_D(|d\phi|^2+ |\psi|^4) <\varepsilon,
 \end{equation}
then
\begin{equation}
\label{5.3} \|\phi\|_{C^k (D_{\frac 1 2})}+\|\psi\|_{C^k (D_{\frac
1 2})} \leq C(\|\n \phi\|_{L^2(D)}+\|\psi\|_{L^4(D)}),
\end{equation}
where $C>0$ is a constant depending only on $k$ and the geometry of
$N$.}

\medskip

 \pr In the sequel, we denote $\|\cdot\|_{L^{k,p}(D)}$ by $|\cdot|_{D,k,p}$.
 If there is no confusion, we may drop the subscript $D$.
 In the proof, $C$ is a constant, varying from line to line
 and $$ D_{\frac 1 2}\subset D^2 \subset D^1 \subset D.$$
  We devide the proof into several steps.
 \medskip
 \par
 \noindent{\bf Step 1.} There is an $\varepsilon>0$ such that
\begin{equation}
\label{t5.3} |d \phi|_{D^1,0,4}\leq C(D^1)
(|d\phi|_{0,2}+|\psi|^2_{0,4}),\quad \forall D^1\subset D,
\end{equation}
where $C(D^1)>0$ is a constant depending only on $D^1$.

 Choose a cut-off function $\eta: 0\leq \eta \leq 1, $ with $\eta|_{D^1}\equiv 1$
 and ${\rm Supp}\eta\subset D$. By (\ref{3.1.11}) we have
\begin{eqnarray*}
|\Delta(\eta\phi)|&\leq&
C(|\phi|+|d\phi|)+|A|_\infty|d\phi|(|d(\eta\phi)|+|\phi d\eta|)
+|\eta\alpha| \\
&\leq&
|A|_\infty|d\phi||d(\eta\phi)|+C(|\phi|+|d\phi|)+|\eta\alpha|,
\end{eqnarray*}
where we denote
$\alpha:=P(A(d\phi(e_\alpha),e_\alpha\cdot\psi);\psi)$. Thus, for
any $p>1$,
\begin{equation}
\label{t5.4} |\Delta(\eta\phi)|_{0,p}\leq
|A|_\infty(|d\phi||d(\eta\phi)|)_{0,p}+C|\phi|_{1,p}+|\eta\alpha|_{0,p}.
\end{equation}
Let $p=\frac{4}{3}$, and without loss of generality we assume
$\int_D\phi=0$ so that $|\phi|_{1,p}\leq C'|d\phi|_{0,p},$ then
$$|A|_\infty(|d\phi||d(\eta\phi)|)_{0,\frac{4}{3}}\leq |A|_\infty |\eta\phi|_{1,4}|d\phi|_{0,2},$$
from this and (\ref{t5.4}) we have
$$
|\eta\phi|_{2,\frac{4}{3}}\leq C(|A|_\infty
|\eta\phi|_{1,4}|d\phi|_{0,2}+|d\phi|_{0,\frac{4}{3}}+|\eta\alpha|_{0,\frac{4}{3}}).
$$
By the Sobolev inequality, $|\eta\phi|_{1,4}\leq
C'|\eta\phi|_{2,\frac{4}{3}}$, so,
\begin{equation}
\label{t5.5} (C^{'-1}-C|A|_\infty
|d\phi|_{0,2})|\eta\phi|_{1,4}\leq
C(|d\phi|_{0,\frac{4}{3}}+|\eta\alpha|_{0,\frac{4}{3}}).
\end{equation}
Moreover,
\begin{eqnarray*}
|\eta\alpha|_{0,\frac{4}{3}}&\leq& C_N(|\psi|^2|\eta d\phi|)_{0,\frac{4}{3}}\\
       &=& C_N(|\psi|^2|d(\eta \phi)-\phi d\eta|)_{0,\frac{4}{3}}\\
&\leq& C|\psi|^2_{0,4}|\eta\phi|_{1,4}+C|\psi|^2_{0,4}.
\end{eqnarray*}
Putting this into (\ref{t5.5}) and choosing $\e$ small, we get
$$|\eta\phi|_{1,4}\leq C(|d\phi|_{0,\frac{4}{3}}+\sqrt{\varepsilon}|\eta\phi|_{1,4}
+|\psi|^2_{0,4}),$$ which yields
$$|\eta\phi|_{1,4}\leq C(|d\phi|_{0,\frac{4}{3}}+|\psi|^2_{0,4})<2\sqrt{\varepsilon}C.$$
\medskip
 \noindent{\bf Step 2.} If $\varepsilon>0$ is small enough,
then
\begin{equation}
\label{t5.6} |\psi|_{D^2,0,q}\leq C(D^2)|\psi|_{D,0,4},\quad
\forall q>1, \quad D^2\subset D^1\subset D,
\end{equation}
where $C(D^2)>0$ is a constant depending only on $D^2$.

Choose a cut-off function $\eta: 0\leq \eta \leq 1, $ with
$\eta|_{D^2}\equiv 1$ and ${\rm
Supp}\eta\subset D^1$. 
Let $\xi=\eta\psi$ have compact support in $D^1$. By the well-known
Lichnerowitz' formula, we have
$$ \partial \hskip -2.2mm \slash^2 \xi=-\Delta
\xi+\frac{1}{4}R\xi=-\Delta \xi,$$ because the scalar curvature
$R\equiv 0$ on $D$. Integrating this yields
\begin{eqnarray*}
|\nabla\xi|_{D^1,0,2} &=& |\partial \hskip -2.2mm \slash \xi|^2_{D^1,0,2} \\
&=& |\partial \hskip -2.2mm \slash (\eta \psi)|^2_{D^1,0,2} =
|\nabla\eta\cdot \psi+\eta\partial \hskip -2.2mm \slash \psi|^2
_{D^1,0,2}\\
&\leq & C(|\psi|^2_{D^1,0,2}+ |\eta\partial \hskip -2.2mm \slash \psi|^2)\\
&\leq & C(|\psi|^2_{D^1,0,2}+ C|d\phi|^2|\eta\psi|^2_{D^1,0,2})\\
 &\leq& C(|\psi|_{D^1,0,2}+|d\phi|^2_{D^1,0,4}|\psi|^2_{D^1,0,4})\\
&\leq& C|\psi|_{D^1,0,4}(1+|d\phi|_{D^1,0,4})\\
&\leq& C'|\psi|_{D^1,0,4}.
\end{eqnarray*}
Hence, Step 2 is proved.
\medskip
\par
\noindent{\bf Step 3.} If $\varepsilon$ is small enough, we have
\begin{equation}
\label{t5.7} |d\phi|_{D^2,0,4}\leq C(D^2)|d\phi|_{D,0,2},\quad
\forall D^2\subset D^1\subset D,
\end{equation}
where $C(D^2)>0$ is a constant depending only on $D^2$.
\medskip
This follows from Steps 1 and 2. For higher order estimates, it
is rather standard. See
for example
\cite{CJLW}. \qed
 \vskip12pt From this theorem we know that a
sequence of Dirac-harmonic maps with small energy has a convergent
subsequence. However, if the energy is large, then a blow-up may
occur. In this case, the concentration of energy may happen.
Namely, $|\n \phi_k(x_k)|^2 \to \infty$ as $k\to \infty$, for a
sequence $x_k\to x_0$ as $k\to \infty$. After a suitable
rescaling, we get a  ``bubble", an entire solution of
(\ref{3.1.11})-(\ref{3.1.10}) with finite energy. By the conformal
invariance of the Dirac-harmonic map, such an entire solution can
be viewed as a Dirac-harmonic map from $\S^2\backslash\{p\}\to N$
with finite energy. In this section, we prove that such a
singularity can be removed as in many conformal problems. Hence,
at the end we obtain a Dirac-harmonic map from $\S^2\to N$. Using the
Theorem 4.3, we can describe the behavior of solutions
$(\phi,\psi)$ near a singular point as  follows: \vskip12pt
\noindent
 {\bf Corollary 4.4.} {\it There is an
$\varepsilon >0$ small enough such that if $(\phi,\psi)$ is a
smooth solution of (\ref{1.4})--(\ref{1.5}) on $D\setminus\{0\}$
with energy $E(\phi,\psi, D)<\varepsilon$, then
 for any $x\in D_{\frac 1 2}$
\begin{equation}
\label{6.2} |d\phi(x)||x|\leq
C(\int_{D(2|x|)}|d\phi|^2)^{\frac{1}{2}},
\end{equation}
\begin{equation}
\label{6.3}
|\psi(x)||x|^{\frac{1}{2}}+|\nabla\psi(x)||x|^{\frac{3}{2}}\leq
C(\int_{D(2|x|)}|\psi|^4)^{\frac{1}{4}},
\end{equation}
  }
   \noindent
   {\it Proof.} Fix any $x_0\in D\setminus\{0\}$, define $\widetilde{\phi}$
   and $\widetilde \psi$ by
   $$\widetilde{\phi}(x):=\phi(x_0+|x_0|x) \quad \hbox{ and }
  \widetilde{\psi}(x):=|x_0|^{\frac{1}{2}}\psi(x_0+|x_0|x). $$
It is clear that $(\widetilde{\phi},\tilde{\psi})$ is a $C^\infty$
solution of (\ref{1.4})--(\ref{1.5}) on $D$ and $E(\widetilde
\phi, \widetilde \psi,D)<\varepsilon$. Applying Theorem 4.3, we
have
$$|d\widetilde{\phi}|_{L^\infty(D_{\frac 1 2})}
\leq C |d\widetilde{\phi}|_{D,0,2}$$ and
\[|\widetilde \psi|_{C^1(D_{\frac 1 2})} \leq C
|d\widetilde{\psi}|_{D,0,4}. \] Scaling back, we prove the
Corollary.\qed
\vskip12pt \noindent {\bf Lemma 4.5.} {\it  Let
$(\phi,\psi)$ be a $C^\infty$ solution of (\ref{1.4})-(\ref{1.5})
on $D\setminus\{0\}$ satisfying $E(\phi,\psi, D)<\varepsilon$,
 then,
\begin{eqnarray}
\label{6.9} \int_0^{2\pi}\frac{1}{r^2}|\phi_\theta|^2d\theta &=&
\int_0^{2\pi}|\phi_r|^2d\theta+\int_0^{2\pi}(\psi,
\partial_r\cdot\psi_r)d\theta  \\
&=&\int_0^{2\pi}|\phi_r|^2d\theta-\int_0^{2\pi}\frac{1}{r^2}(\psi,
\partial_\theta\cdot\psi_\theta)d\theta,\nonumber
\end{eqnarray}
where $(r,\theta)$ are the polar coordinates in $D$ centered at 0,
$\psi_r:=\widetilde{\nabla}_{\partial_r}\psi$.
 }
\vskip12pt \noindent {\it Proof.} From Proposition 3.2,
\begin{eqnarray*}
T&=&[|\phi_x|^2-|\phi_y|^2-2i \la \phi_x,\phi_y \ra ]+[ \la
\psi,\partial_x\cdot \widetilde{\nabla}_{\partial_x}\psi\ra -i \la
\psi,\partial_x\cdot
\widetilde{\nabla}_{\partial_y}\psi) \ra ]\\
&=:& (A-iB)+(A'-iB'), \quad z=x+iy\in D
\end{eqnarray*}
is holomorphic on $D\setminus\{0\}$. 
  By Corollary 4.4, we know that
$$|A-iB|\leq 2(|\phi_x|^2+|\phi_y|^2)\leq C|z|^{-2}.$$
Noting that
$\wa\psi=(\nabla_{e_\alpha}\psi^i+\Gamma^i_{jk}\phi^j_\alpha\psi^k)\otimes\partial_{y^i}$,
we have
$$|\widetilde{\nabla}\psi|\leq C(|\nabla\psi|+|d\phi||\psi|).$$
 By  Corollary 4.4 again, we have
$$|A'-iB'|\leq 2|\psi||\widetilde{\nabla}\psi|\leq
C(|\psi||\nabla\psi|+|d\phi||\psi|^2)\leq C|z|^{-2}.$$ Therefore
$|T(z)|\leq C|z|^{-2}.$ Furthermore,
$$ \int_D|A'-iB'|\leq 2\int_D|\psi||\widetilde{\nabla}\psi|<\infty$$
and  $\int_D|A-iB|\leq 2\int_D|d\phi|^2<\infty.$ Thus,
$\int_D|T(z)|<\infty,$ which implies that $T(z)$ has a pole at
$z=0$ of order at most one.
 Hence, $zT(z)$ is holomorphic on $D$ and
\begin{equation}
\label{6.11} 0={\rm Im}[\int_{|z|=r}zT(z)dz]=\int_0^{2\pi}{\rm
Re}[z^2T(z)]d\theta.
\end{equation}
It is easy to compute that
\begin{eqnarray}
\label{6.12} {\rm Re}[z^2T(z)]&=&
r^2[(A\cos 2\theta+B\sin 2\theta)+(A'\cos 2\theta+B'\sin 2\theta)]\\
&=&\nonumber r^2|\phi_r|^2-|\phi_\theta|^2- \la \psi,
\partial_\theta\cdot \psi_\theta\ra.
\end{eqnarray}
 Now the desired equalities follow. \hfill $\Box$ \vskip12pt \noindent
{\bf Remark 4.6.} Integrating (\ref{6.9}) yields:
\begin{equation}
\label{6.15}
\int_D|\phi_r|^2-\int_D\frac{1}{r^2}|\phi_\theta|^2=-\int_D \la
\psi,\partial_r \cdot \psi_r \ra :=I.
\end{equation}
 From (\ref{6.9}) and
$$\int_{|z|=r}|\phi_r|^2+\int_{|z|=r}\frac{1}{r^2}|\phi_\theta|^2=\int_{|z|=r}|d\phi|^2:=E_r(\phi)
$$
we have
\begin{equation}
\label{6.16} \int_{|z|=r}|\phi_r|^2=\frac{1}{2}E_r+\frac{1}{2}I_r,
\end{equation}
and
\begin{equation}
\label{6.17}
\int_{|z|=r}\frac{1}{r^2}|\phi_\theta|^2=\frac{1}{2}E_r-\frac{1}{2}I_r,
\end{equation}
where $I_r:=-\int_{|z|=r}(\psi,\partial_r\cdot \psi_r).$

 Now we
can state the following
 \vskip12pt
\noindent {\bf Theorem 4.7} ({\bf Removable singularity theorem})
{\it Let $(\phi,\psi)$ be a solution of (\ref{1.4}) and
(\ref{1.5}) which is $C^\infty$ on $U\setminus\{p\}$ for some
$p\in U\subset M$. If $(\phi,\psi)$ has finite energy, then
$(\phi,\psi)$ extends to a $C^\infty$ solution on $U$. }
\vskip12pt \noindent {\it Proof.} By rescaling, we may assume that
\begin{equation}
\label{6.22} \int_{D(2)}(|d\phi|^2+|\psi|^4)<\varepsilon.
\end{equation}
Choose a function $q(r)$ on $D$ which is piecewise linear in
$\log r$ with
$$q(2^{-m})=\frac{1}{2\pi}\int_0^{2\pi}\phi(2^{-m},\theta)d\theta,$$
then  we have (cf. \cite{SU})
\begin{equation}
\label{6.25}
\int_D|dq-d\phi|^2=\int_{r=1}(q-\phi)\phi_r-\int_D(q-\phi)\Delta
(q-\phi)
\end{equation}
with
\begin{equation}
\label{6.26}
|q-\phi|_\infty:=|q-\phi|_{C^0(D)}<2^3\sqrt{\varepsilon}.
\end{equation}
Using (3.21), we have $\Delta (q-\phi)=-\Delta
\phi=-A(\phi)(d\phi,d\phi)-\alpha,$ so,
\begin{eqnarray}
\label{6.27}
|\int_D(q-\phi)\Delta (q-\phi)|&\leq & |q-\phi|_\infty|A|_\infty\int_D|d\phi|^2+|q-\phi|_\infty\int_D|\alpha|\nonumber \\
&<&
\delta\int_D|d\phi|^2+C\sqrt{\varepsilon}\int_D|\psi|^2|d\phi|,
\end{eqnarray}
where $2^3\sqrt{\varepsilon}|A|_\infty<\delta$ for some constant
$\delta >0$ small, and $C>0$ is a constant. As for the first term
on the RHS of (\ref{6.25}),
 \begin{eqnarray}
\label{6.28}
\int_{r=1}(q-\phi)\phi_r &\leq & (\int_{r=1}|q-\phi|^2)^{\frac{1}{2}}(\int_{r=1}|\phi_r|^2)^{\frac{1}{2}}\nonumber \\
&\leq &
(\int_{r=1}|\phi_\theta|^2)^{\frac{1}{2}}(\int_{r=1}|\phi_r|^2)^{\frac{1}{2}}.
\end{eqnarray}
On the other hand,
by (4.16),
\begin{equation}
\label{6.29} \int_D|dq-d\phi|^2\geq
\int_D\frac{1}{r^2}|\phi_\theta|^2=\frac{1}{2}E(\phi)-\frac{1}{2}I,
\end{equation}
where $I:=-\int_D \la \psi,\partial_r\cdot \psi_r\ra.$ Inserting
(\ref{6.27}), (\ref{6.28}) and (\ref{6.29}) into (\ref{6.25}), we
get
\begin{eqnarray}
\label{6.30} \frac{1}{2}\int_D|d\phi|^2-\frac{1}{2}I&\leq
&(\int_{r=1}|\phi_\theta|^2)^{\frac{1}{2}}(\int_{r=1}|\phi_r|^2)^{\frac{1}{2}}
\nonumber \\
&&+\delta\int_D|d\phi|^2+C\sqrt{\varepsilon}\int_D|\psi|^2|d\phi|.
\end{eqnarray}
By (\ref{6.16}) and (\ref{6.17}), we have
\begin{eqnarray*}
(\int_{r=1}|\phi_\theta|^2)^{\frac{1}{2}}(\int_{r=1}|\phi_r|^2)^{\frac{1}{2}}
&=& (\frac{1}{2}E_1-\frac{1}{2}I_1)^\frac{1}{2}(\frac{1}{2}E_1+\frac{1}{2}I_1)^\frac{1}{2}\nonumber \\
&\leq& \frac{1}{2}E_1=\frac{1}{2}\int_{r=1}|d\phi|^2.
\end{eqnarray*}
It follows
\begin{equation}
\label{6.31} (1-2\delta)\int_D|d\phi|^2\leq \int_{r=1}|d\phi|^2+
[2C\sqrt{\varepsilon}\int_D|\psi|^2|d\phi|- \int_D \la
\psi,\partial_r\cdot \psi_r \ra ].
\end{equation}
By a scaling argument, this yields
\begin{eqnarray}
\label{6.32}
 (1-2\delta)\int_{D_r}|d\phi|^2&\leq & r\int_{\partial D_r}|d\phi|^2+[2C\sqrt{\varepsilon}\int_{D_r}|\psi|^2|d\phi|-
\int_{D_r} \la \psi,\partial_r\cdot \psi_r \ra ] \nonumber \\
&\leq & r\int_{\partial
D_r}|d\phi|^2+C\sqrt{\varepsilon}\int_{D_r}|\psi|^4+C\sqrt{\varepsilon}\int_{D_r}|d\phi|^2
+\int_{D_r}|\psi||\nabla\psi|\nonumber \\
&\leq&r\int_{\partial
D_r}|d\phi|^2+C\sqrt{\varepsilon}\int_{D_r}|d\phi|^2+C\int_{D_r}|\psi|^4+C\int_{D_r}|\nabla
\psi|^{\frac 4 3}.
\end{eqnarray}
Now we need the following lemma. It is the elliptic estimate with
boundary. For completeness we give a proof here.
\par
\vskip12pt \noindent {\bf Lemma 4.8.}  {\it Let $u$ be a complex
function satisfying
\begin{equation}
\label{6.33.1} \left\{
\begin{array}{l} \bar \partial u=f,
 \quad {\rm in}\quad D,
\\ u|_{\partial D}=\varphi,
\end{array}
\right.
\end{equation}
with $\varphi \in L^{1,p}(\partial D)$ and $f\in L^p(D)$ for some
$p>1$, where $D$ is the unit disc on $\R^2$ centered at the
origin, then the following estimate holds
\begin{equation}
\label{6.34.1} |u|_{D,1,p}\leq C(|f|_{D,0,p}+|\varphi|_{\partial
D, 1,p}).
\end{equation}
If instead $u$ satisfies
\begin{equation}
\label{6.35.1} \left\{
\begin{array}{l}  \partial u=f,
 \quad {\rm in}\quad D,
\\ u|_{\partial D}=\varphi,
\end{array}
\right.
\end{equation}
 then the same estimate holds.}
\vskip12pt \noindent {\it Proof.} We first consider the following
boundary value problem:
\begin{equation}
\label{6.52} \left\{
\begin{array}{l} \bar \partial w_1=0,\\{\rm Re}w_1|_{\partial D}=\varphi, \quad \int_{\partial D}{\rm Im}w_1=2\pi
c_0,
\end{array}
\right.
\end{equation}
where $c_0:=\frac {1}{ 2\pi} \int_{\partial D} {\rm Im} \varphi.$
 It is clear that (cf. \cite{Be}, Theorem 38):
\begin{equation}
\label{6.53} |w_1|_{C^\alpha(D)}+|\nabla w_1|_{L^p(D)}\leq C
(|\varphi|_{\partial D,1,p}+c_0),
\end{equation}
where $C>0$ is a constant.
\par
Next, let $w_2:=u-w_1,$ then it satisfies
\begin{equation}
\label{6.55} \left\{
\begin{array}{l} \bar \partial w_2=f \quad {\rm in}\quad D,
\\{\rm Re}w_2|_{\partial D}=0, \quad \int_{\partial D}{\rm
Im}w_2=0.
\end{array}
\right.
\end{equation}
The elliptic estimates yield that
\begin{equation}
\label{6.56} |w_2|_{L^{1,p}(D)}\leq C |f|_{L^p(D)}.
\end{equation} In fact, by the Schwarz-Poisson formula (see Theorem 21 in
\cite{Be}), one has
\begin{equation}
\label{6.59} w_2=-\frac 1 \pi \int_{|\zeta|<1}(\frac {f}{\zeta
-z}+\frac{\bar z\bar f}{1-z\bar \zeta})d\xi d\eta.
\end{equation}
 From this and the boundedness of the Riesz transformation, we
obtain (\ref{6.56}). Combining (\ref{6.53}) and (\ref{6.56}) then
gives the desired estimate (\ref{6.34.1}). The proof is similar if
$u$ satisfies (\ref{6.35.1}). \hfill  $\Box$
\par
Now we return to the proof of Theorem 4.7. Recall (3.22)
\begin{equation}
\label{6.59.1}
\partial \hskip -2.2mm \slash
 \psi = A(d\phi(e_\a), e_\a\cdot \psi )\quad {\rm in}\quad D\setminus \{0\}.
\end{equation}
We choose a cut-off function $\eta_\varepsilon\in
C^\infty_0(D_{2\varepsilon})$ such that $\eta_\varepsilon=1$ in
$D_\varepsilon(0)$ and $|d\eta_\varepsilon|<C/\varepsilon.$ Then
we have
$$ \partial \hskip -2.2mm \slash ((1-\eta_\varepsilon)\psi)=(1-\eta_\varepsilon) A(d\phi(e_\a), e_\a\cdot \psi )
-\nabla\eta_\varepsilon\cdot\psi.  $$
 From Lemma 4.8, we have
\begin{equation}
\label{6.60} |(1-\eta_\varepsilon)\psi|_{D,1,\frac 4 3}\leq
C|d\phi|_{D,0,2}|\psi|_{D,0,4}+C|\psi|_{\partial D,1,\frac 4
3}+C|\nabla\eta_\varepsilon\cdot\psi|_{D,0,\frac 4 3}.
\end{equation}
Letting $\varepsilon\to 0$, using
$$\lim _{\varepsilon\to 0}\frac 1 {\varepsilon
^{\frac 43}}
\int_{D_{2\varepsilon}}|\psi|^{\frac 4 3}=0,$$
 the smallness of
$|d\phi|_{D,0,2}$ and  the Sobolev embedding theorem, we obtain
$$(\int_D|\psi|^4)^{\frac 1 4}\leq C (\int_{\partial
D}|\nabla\psi|^{\frac 4 3})^{\frac 3 4}+C(\int_{\partial
D}|\psi|^4)^{\frac 1 4}.$$ By rescaling, we have for any $0\leq
r\leq 1$
\begin{eqnarray*}
(\int_{D_r}|\psi|^4)^{\frac 1 4}&\leq & C (r\int_{\partial
D_r}|\nabla\psi|^{\frac 4 3})^{\frac 3 4}+C(r\int_{\partial
D_r}|\psi|^4)^{\frac 1 4}\\
&\leq& C (r\int_{\partial D_r}|\nabla\psi|^{\frac 4 3})^{\frac 1
4}+C(r\int_{\partial D_r}|\psi|^4)^{\frac 1 4}.
\end{eqnarray*}
Thus,
\begin{equation}
\label{6.61} \int_{D_r}|\psi|^4 \leq C r\int_{\partial
D_r}|\nabla\psi|^{\frac 4 3}+Cr\int_{\partial D_r}|\psi|^4.
\end{equation}
Let $\bar{\psi}:=\frac{1}{2\pi}\int_{\partial D}\psi.$ Note that
$$\partial \hskip -2.2mm \slash
 (\psi-\bar{\psi})=A(d\phi(e_\a), e_\a\cdot (\psi-\bar{\psi}) )+A(d\phi(e_\a), e_\a\cdot
 \bar{\psi})\quad {\rm in}\quad D\setminus\{0\}.
 $$
 By an argument similar to the one used in obtaining (\ref{6.60}) and
 using  the Poincar\'e inequality,
 we have
 \begin{eqnarray}
 \label{6.62}
 |\psi-\bar{\psi}|_{D,1,\frac 4 3}&\leq& C|d\phi|_{D,0,2}|\psi-\bar{\psi}|_{D,1,\frac 4
 3}+C(|d\phi||\bar \psi|)_{D,0,\frac 4 3}\nonumber \\
 &&+C|\nabla\psi|_{\partial D,0,\frac 4 3}.
 \end{eqnarray}
 Again, by the smallness of $|d\phi|_{D,0,2}$ we obtain
\begin{eqnarray*}
(\int_D|\nabla\psi|^{\frac 4 3})^{\frac 3 4}&\leq &
C(\int_{\partial D}|\nabla \psi|^{\frac 4 3})^{\frac 3
4}+C|\bar{\psi}|(\int_D|d\phi|^2)^{\frac 1 2}\\
&\leq& C(\int_{\partial D}|\nabla \psi|^{\frac 4 3})^{\frac 3
4}+C(\int_{\partial D}|\psi|^4)^{\frac 1 4}(\int_
D|d\phi|^2)^{\frac 1 2}.
\end{eqnarray*}
So we have
\begin{eqnarray*}
\int_D|\nabla\psi|^{\frac 4 3}&\leq & C\int_{\partial D}|\nabla
\psi|^{\frac 4 3}+C(\int_{\partial D}|\psi|^4)^{\frac 1
3}(\int_D|d\phi|^2)^{\frac 2 3}\\
&\leq& C\int_{\partial D}|\nabla \psi|^{\frac 4
3}+\varepsilon_1\int_D|d\phi|^2+\frac{C}{\varepsilon_1}\int_{\partial
D}|\psi|^4,
\end{eqnarray*}
where $\varepsilon_1>0$ is a small constant. Hence, for $0\leq r
\leq 1$,
\begin{equation}
\label{6.63} \int_{D_r}|\nabla\psi|^{\frac 4 3}\leq
Cr\int_{\partial D_r}|\nabla \psi|^{\frac 4
3}+\varepsilon_1\int_{D_r}|d\phi|^2+\frac{Cr}{\varepsilon_1}\int_{\partial
D_r}|\psi|^4.
\end{equation}
Putting (\ref{6.32}), (\ref{6.61}) and (\ref{6.63}) together, we
have for any $0\leq r \leq 1$ and some constant $C>0$
\begin{equation}
\label{6.64}
\int_{D_r}|d\phi|^2+\int_{D_r}|\psi|^4+\int_{D_r}|\nabla\psi|^{\frac
4 3}\leq Cr(\int_{\partial D_r}|d\phi|^2+\int_{\partial
D_r}|\psi|^4+\int_{\partial D_r}|\nabla\psi|^{\frac 4 3}).
\end{equation}
Denote
$F(r):=\int_{D_r}|d\phi|^2+\int_{D_r}|\psi|^4+\int_{D_r}|\nabla\psi|^{\frac
4 3}$. (\ref{6.64}) implies that
\begin{equation}
\label{6.65} F(r)\leq C r F'(r).
\end{equation}
Integrating this inequality yields
\begin{equation}
\label{6.66} F(r)\leq  F(1)r^{\frac 1 C}.
\end{equation}
 From this we can easily conclude that there are some $\beta>1$ and
$2>q>\frac 4 3$ such that
\begin{equation}
\label{b.67} \phi\in L^{1,2\beta}(D),\quad \psi\in L^{1,q}(D).
\end{equation}

\par Now let us consider $\psi$. Recall that $\partial \hskip-2.2mm
\slash \psi = A(d\phi(e_\a), e_\a\cdot \psi ):=a,$ and note that
\begin{eqnarray*}
\int_{D_{\frac 1 2}}|a|^{q_1} &\leq & C\int_{D_{\frac 1 2}}(|d\phi||\psi|)^{q_1}\\
&\leq & C(\int_{D_{\frac 1
2}}|d\phi|^{q_1p'})^{\frac{1}{p'}}(\int_{D_{\frac 1
2}}|\psi|^{q_1q'})^{\frac{1}{q'}},
\end{eqnarray*}
where $q_1>0, p'$ and $q'$ are constants chosen as follows
$$q_1=\frac{2\beta q}{2\beta-(\beta-1)q},\quad p'=\frac{2q}{2q-(2-q)q_1}, \quad q'=\frac{2q}{(2-q)q_1}.$$
Since $\beta>1$, we have $q_1>q$ and
\begin{equation}
\label{6.38}
 \int_{D_{\frac 1 2}}|a|^{q_1}\leq C(\int_{D_{\frac 1 2}}|d\phi|^{2\beta}
)^{\frac{1}{p'}}(\int_{D_{\frac 1
2}}|\psi|^{q^*})^{\frac{1}{q'}}<\infty,
\end{equation}
where $q^*:=\frac{2q}{2-q}$. We note that $\psi\in L^{q_1}$
because of $q_1<q^*$.

 From the regularity theory of Cauchy-Riemann operators
we have
\begin{equation}
\label{6.39} \psi\in L^{q_1}_1 \quad {\rm with} \quad
q_1=\frac{2\beta q}{2\beta-(\beta-1)q}>q.
\end{equation}

 From above, we see that $\phi\in L^{1,2\beta}$
and $\psi\in L^{1,q} (4/3<q<2)$ imply $\phi\in L^{1,2\beta}$ and
$\psi\in L^{1,q_1}.$ By iteration, $\phi\in L^{1,2\beta}$ and
$\psi\in L^{1,q_n}$ imply $\phi\in L^{1,2\beta}$ and $\psi\in
L^{1,q_{n+1}}$ with
$$q_{n+1}=\frac{2\beta q_n}{2\beta-(\beta-1) q_n}.$$
Since
$$\frac{q_{n+1}}{q_n}=\frac{2\beta }{2\beta-(\beta-1)q_n}>\frac{3\beta }{\beta+2}>1,$$
there exists $q_n$ such that $q_n> 2.$ Therefore, we have that
$\phi\in L^{1,p_0}(D)$ and $\psi\in L^{1,q_0}(D)$ for some $p_0>2$
and $q_0>2$. We can then conclude that $(\phi,\psi)$ is smooth on
$D$ through the standard bootstrap method. We omit the details
here. This completes the proof of the theorem. \hfill $\Box$

\end{document}